\documentclass[twocolumn,amsmath,showpacs,amssymb,aps,secnumarabic,%
    nofootinbib,superscriptaddress]{revtex4}
\usepackage{times,mathptmx,amsthm,graphicx}

\newtheorem{theorem}{Theorem}
\newtheorem{corollary}[theorem]{Corollary}
\theoremstyle{definition}
\newtheorem*{definition}{Definition}

\newcommand{\R}{\mathbb{R}}
\newcommand{\Z}{\mathbb{Z}}
\newcommand{\cC}{\mathcal{C}}
\newcommand{\cB}{\mathcal{B}}
\newcommand{\cS}{\mathcal{S}}
\newcommand{\wcB}{\widetilde{\cB}}
\newcommand{\wcC}{\widetilde{\cC}}
\newcommand{\wcS}{\widetilde{\cS}}
\newcommand{\vv}{\vec{v}}
\newcommand{\ie}{\emph{i.e.}}

\makeatletter
\def\@pacs@name{MSC numbers: }%
\makeatother

\begin{document}

\title[Lattice packings with gap defects]
{Lattice packings with gap defects are not completely saturated}

\author{Greg~Kuperberg}
\email{greg@math.ucdavis.edu}
\affiliation{Department of Mathematics, University of California, Davis, CA 95616}

\author{Krystyna~Kuperberg}
\email{kuperkm@math.auburn.edu}
\affiliation{Department of Mathematics, Auburn University, Auburn, AL 36849}

\author{W{\l}odzimierz~Kuperberg}
\email{kuperwl@math.auburn.edu}
\affiliation{Department of Mathematics, Auburn University, Auburn, AL 36849}



\thanks{G.~Kuperberg's research was supported in part by the National
Science Foundation grant \#DMS-0072342}

\thanks{K.~Kuperberg's research was supported in part by the National
Science Foundation grant \#DMS-0204081}

\begin{abstract} We show that a honeycomb circle packing in $\R^2$ with a
linear gap defect cannot be completely saturated, no matter how
narrow the gap is.  The result is motivated by an open problem of G.~Fejes
T\'oth, G.~Kuperberg, and W.~Kuperberg, which asks whether of a honeycomb
circle packing with a linear shift defect is completely saturated. We also show
that an fcc sphere packing in $\R^3$ with a planar gap defect is also not
completely saturated.
\end{abstract}

\maketitle

\section{Introduction} A packing of unit spheres in $\R^n$ is
\emph{$k$-saturated} if it is not possible to replace $k-1$ spheres by $k$ and
still have a packing; it is \emph{completely saturated} if it is $k$-saturated
for all $k$ \cite{FTK}. As discussed in \cite{FTK}, every completely saturated
packing has maximum density, and every periodic packing with maximum density is
completely saturated.  Thus the densest lattice packing of circles in $\R^2$
(the \emph{honeycomb} circle packing) is completely saturated.  This
motivates the question:  Is the honeycomb circle packing the only
one which is completely saturated?

\begin{figure}[h]
\begin{center}
\includegraphics{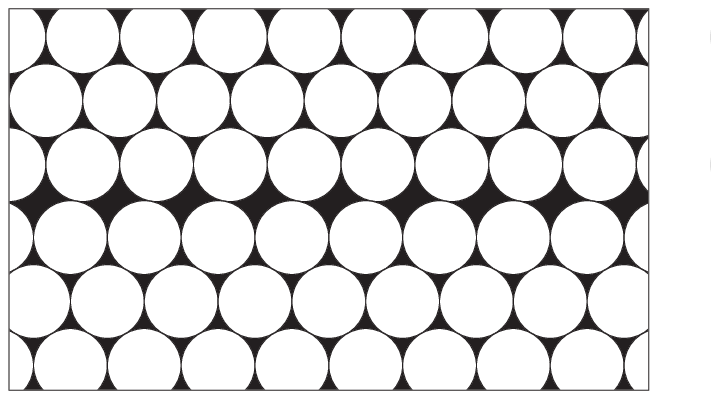}
\end{center}
\caption{Honeycomb circle packing with a shift defect.}
\label{shift}
\end{figure}

G.~Fejes T\'oth, G.~Kuperberg, and W.~Kuperberg \cite{FTKK} asked whether a
honeycomb circle packing with a linear shift defect, as in Figure~\ref{shift},
is completely saturated. This question appears to be closely related to the
(still open) conjecture of L.~Fejes T\'oth \cite{LFT} asserting that the same
packing is \emph{solid}, meaning that if you obtain another packing by
rearranging finitely many circles, it is congruent to the original one.

In this paper we will instead consider packings with a linear or
planar gap defect in which the spheres on either side do not touch.

\begin{definition}\label{n-gap} Let $\cS$ be a sphere packing in $\R^n$ let $H
\subset \R^n$ be a hyperplane, let $H^+$ be a closed half-space bounded by $H$,
and let $\vv$ be a vector perpendicular to $H$ in the direction of $H^+$. Let
$\wcS$ be the packing obtained from $\cS$ by moving all spheres with centers in
$H^+$ by the vector $\vv$.  (Note that $\wcS$ is a packing because the motion
does not decrease the distances between sphere centers.) Then $\wcS$ is a \emph{sphere
packing with a hyperplane gap defect} and $\|\vv\|$ is the \emph{width} of
the gap.
\end{definition} 

\begin{figure}[h]
\begin{center}
\includegraphics{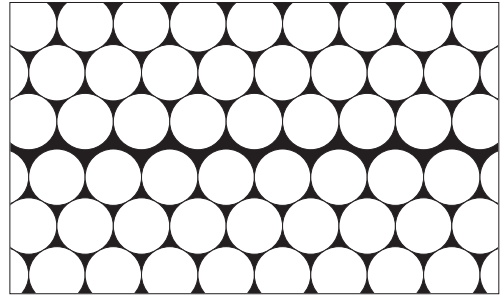}
\end{center}
\caption{Honeycomb circle packing with a gap defect.}
\label{gapa}
\end{figure}

Figure~\ref{gapa} shows a honeycomb circle packing with a linear gap defect.

\begin{theorem} A honeycomb circle packing in $\R^2$ with a linear gap defect
is not completely saturated, regardless of the width and direction of the gap. 
An fcc (face-centered cubic lattice) packing of spheres in $\R^3$ with an
arbitrary planar gap defect is also not completely saturated.
\end{theorem}

Since the honeycomb circle packing in $\R^2$ \cite{LA} and the fcc sphere
packing in $\R^3$ \cite{GA} are the unique densest lattice sphere packings in
2 and 3 dimensions, we obtain the following corollary.

\begin{corollary} In $\R^2$ and $\R^3$, no lattice sphere packing with a
hyperplane gap defect is completely saturated.
\end{corollary}

\section{Honeycomb packing with gap defects}

\subsection{Gap along a lattice line}
\label{latgap}

Let $\cC$ be the honeycomb circle packing in $\R^2$ with a circle center at
every point of the form $(2i+j,\sqrt{3}j)$, with $i,j \in \Z$. Let $d > 0$ and
move all circles with center on or above the $x$-axis by the vector $\vv =
(0,d)$ to obtain the packing $\wcC$. 

We will show that $\wcC$ is not completely saturated.  It suffices to rearrange
finitely many circles in $\wcC$ to create a void with enough space for an extra
circle.  To achieve that, we will widen the gap by a factor bounded away from
$1$ in a sufficiently large finite region.  If we can create a gap of width $d'
= d+\delta$, with $\delta$ a non-decreasing function of $d$, then by repeating
the operation we can create a gap of any width.  In particular, if the 
width is at least 2, there is room for another circle.

\begin{figure}[ht]
\begin{center}
\scalebox{.85}{\includegraphics{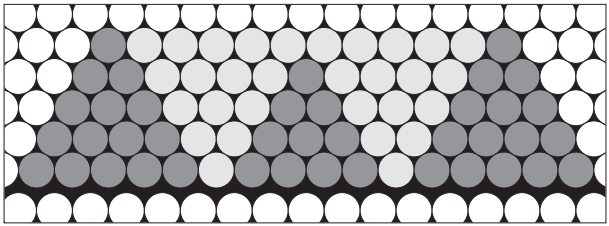}}
\end{center}
\caption{Five triangular blocks.}
\label{blocks}
\end{figure}

\begin{figure}[ht]
\begin{center}
\scalebox{.85}{\includegraphics{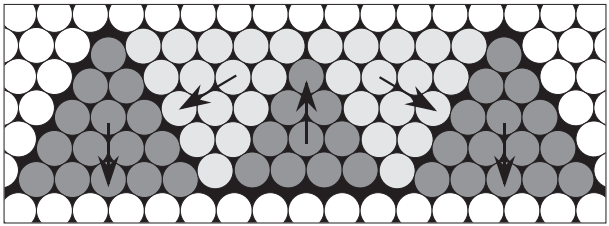}}
\end{center}
\caption{The rearrangement in three moves.}
\label{move}
\end{figure}

We define a \emph{triangular block of size $n$} to be a set of \break
$n(n+1)/2$ circles whose centers lie in an equilateral triangle of edge length
$2n-2$, \ie, with $n$ circles along each edge. Note that a triangular block in
$\wcC$ must lie entirely on one side of the gap.  If the block points up, we
will call it a $\Delta_n$-block; otherwise we will call it a $\nabla_n$-block. 
The rearrangement to widen the gap is as follows:  Take two $\Delta_n$-blocks,
two $\nabla_n$-blocks, and one $\Delta_{n-1}$-block that form a trapezoid, as
in Figure~\ref{blocks}.  Move the two outer $\Delta_n$-blocks down by a
distance of $d_1 \le d$, then move the two $\nabla_n$-blocks at $60$ degree
angles from vertical by $d_2$, then move the middle $\Delta_n$-block up by
$d_3$, as in Figure~\ref{move}.   If we choose the distances $d_1$, $d_2$, and
$d_3$ to maximize $d_3$, then $d_3$ is at least weakly monotonic in $d$,
because increasing $d$ relaxes the constraints on the other parameters. Thus
the gap is widened by $\delta$ monotonic in $d$, as desired.  Since the
parameter $n$ is arbitrary, it can be taken large enough to repeat the
operation with smaller blocks to reach any desired width.

\subsection{Gap along an arbitrary line}
\label{linegap}

\begin{figure}[ht]
\begin{center}
\includegraphics{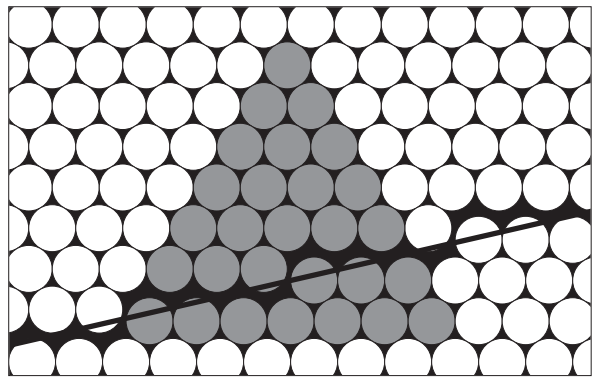}
\end{center}
\caption{A triangular block split by the gap.}
\label{split}
\end{figure}

Let $\cC$ be the above honeycomb circle packing and let $\ell$ be any line, and
now let $\wcC$ be the packing $\cC$ with a gap along $\ell$. Our plan in this
case is to create another (finite) gap parallel to a lattice line of $\cC$ and
reduce to the case of Section~\ref{latgap}.  The gap can be made any length
without sacrificing width.  For each $n \ge 3$, we can choose a triangular
block $T$ in $\cC$ of length $n$ so that $\ell$ intersects the triangle formed
by its circle centers, moreover so that $\ell$ does not meet the vertices. 
Then $\ell$ divides $T$ into two sub-blocks $T_1$ and $T_2$, one of which, say
$T_1$, has only one corner of $T$. In $\wcC$, these two sub-blocks are
separated by the gap, as in Figure~\ref{split}. In $\wcC$, then, $T_1$ can be
moved towards $T_2$ to introduce gaps on its other two sides. The length $n$ is
arbitrary and the width of the new gaps does not depend on $n$.  This completes
the reduction.

\section{The fcc packing with a gap}

The fcc sphere packing in $\R^3$ decomposes into honeycomb layers, and, as
noted by Kepler \cite{JK}, it also decomposes into square-lattice layers.  We
prove below that the fcc-lattice sphere packing with a gap defect is never completely
saturated, no matter how narrow the gap is and no matter along which plane the
gap is formed. As before, we first establish special cases.

\subsection{Gap parallel to a square layer}
\label{squaregap}

Let $\cB$ be the fcc packing in $\R^3$ whose sphere centers are at
$(2i+k,2j+k,\sqrt{2}k)$ with $i,j,k \in \Z$. Let $\wcB$ be the packing $\cB$
with sphere centers in the upper half-space moved by $\vv = (0,0,d)$ for some
$d>0$.  In this model the square-lattice layers are parallel to coordinate
planes of $\R^3$, so $\wcB$ has a square-layer gap.

\begin{figure}[ht]
\begin{center}
\scalebox{.75}{\includegraphics{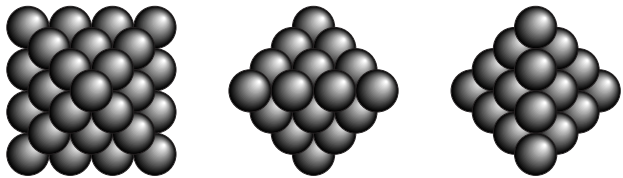}}
\end{center}
\caption{A square pyramid block and two simplex blocks (top view).}
\label{clusters}
\end{figure}

To show that $\wcB$ is not completely saturated, we will delete finitely many
spheres, then rearrange finitely many others to create a void large enough to
accommodate more spheres than the ones deleted.  We define a \emph{square
pyramid block of size $n$} to be a collection of $n(n+1)(2n+1)/6$ spheres in
$\wcB$ whose sphere centers lie in a square pyramid of edge length \break
$2n-2$, \ie, with $n$ spheres along each of its edges. We define  a
\emph{simplex block of size $n$} to be a collection of $n(n+1)(n+2)/6$ spheres
in $\wcB$ whose spheres centers lie in a regular simplex of edge length $2n-2$.
Both kinds of blocks have $n$ spheres along each edge; see
Figure~\ref{clusters} for examples.

\begin{figure}[ht]
\begin{center}
\includegraphics{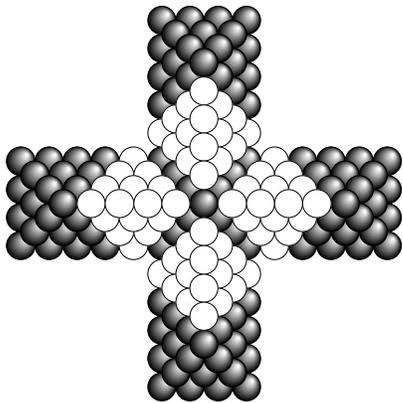}
\end{center}
\caption{Nine blocks forming a cross gable roof.}
\label{config}
\end{figure}

By analogy with Section~\ref{latgap}, we will widen the gap by moving 5 square
pyramid blocks and 4 simplex blocks forming a cross gable roof, as in
Figure~\ref{config}.  The bases of the square pyramid blocks lie on the gap. We
move the four outer square pyramid blocks down, then move the simplex blocks
down and out, then move the middle square pyramid up after removing its top
sphere.  As before, we choose the motions in order to maximize the
size of the last motion; in this arrangement the last motion is weakly
monotonic in $d$.  Thus the gap widens by $\delta$ with $\delta$ weakly
monotonic in $d$.  As before, $\delta$ does not depend on $n$.

As in Section~\ref{latgap}, $n$ can be large enough to iterate
the procedure many times.  Let $k$ be the number of iterations needed
to widen the gap to 2; $k$ only depends on $d$.  One sphere is 
deleted in each iteration, for a total of $k$ spheres removed.
Since the void at the end grows with $n$, we can choose $n$ large
enough to accommodate these $k$ spheres and then some.

\subsection{Gap parallel to a honeycomb layer}

\begin{figure}[ht]
\begin{center}
\includegraphics{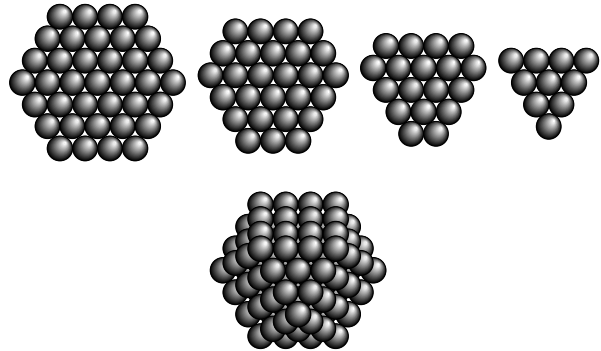}
\end{center}
\caption{A half-cuboctahedron block, exploded into layers and assembled.}
\label{hex}
\end{figure}

Our strategy in this case is to expose a square-layer gap in a large finite
region by moving a large block that abuts the honeycomb gap.  This reduces the
problem to the case considered in Section~\ref{squaregap}. The shape of the block
is a cuboctahedron sawed in half to expose a regular hexagon face, as in
Figure~\ref{hex}.  As the figure indicates, the three square faces lie in
square-lattice directions.  Moving the block creates gaps along these faces.

\subsection{Gap parallel to an arbitrary plane}

Assume now that the plane of the gap does not belong to any of the four
parallel classes of the honeycomb layers. In a manner similar to that of
Section \ref{linegap}, we construct a large tetrahedron block which is cut by
the gap into two sub-blocks $T_1$ and $T_2$.  If the plane is not  parallel to
a square-lattice layer, then we can arrange that $T_1$ contains 1 vertex of the
tetrahedron and $T_2$ contains 3 vertices. The sub-block $T_1$ can then be
moved towards the opposite face, which creates a honeycomb-layer gap and
reduces the problem to that of Section \ref{squaregap}.

\section{Open problems and conjectures}

\begin{enumerate}
\item  Is there a completely saturated sphere packing with a gap defect in
$\R^n$ ($n\geq 2$)?

\item Let $\cS$ be the densest lattice packing of spheres in $\R^4$, which is
composed of fcc layers.  Let $H$ be a hyperplane in an fcc direction and let
$\widehat{\cS}$ be a sphere packing obtained from $\cS$ by introducing a gap
defect along $H$. We conjecture that $\widehat{\cS}$ is not completely saturated.

\item  Let $\mathcal{T}$ be a tiling of $\R^2$ with regular hexagons, 
and let $\widehat{\mathcal{T}}$ be a packing obtained from $\mathcal{T}$
by introducing a gap.  If the gap is sufficiently narrow, is
$\widehat{\mathcal{T}}$ completely saturated?

\end{enumerate}


\begin{thebibliography}{10}

\bibitem{FTKK} Fejes T{\'o}th,~G{\'a}bor, Kuperberg,~Greg, and
Kuperberg,~W{\l}o\-dzi\-mierz, \emph{Highly saturated packings and reduced
coverings\/}, Monatsh. Math. 125(1998), no. 2, 127-145.

\bibitem{FTK} Fejes T{\'o}th,~G{\'a}bor and
Kuperberg,~W{\l}odzimierz, \emph{Packing and covering with convex
sets}, Chapter 3.3 in: Handbook of Convex Geometry (P.M.~Gruber and
J.M.~Wills, Eds.), Elsevier 1993, 799-860.

\bibitem{LFT} Fejes T{\'o}th, L{\'a}szl{\'o}, \emph{Solid
circle-packings and circle-coverings}, Studia Sci.~Math.~Hungar.
3(1968), 401-409.

\bibitem{GA} Gauss, Carl Friedrich, \emph{Untersuchungen \"uber die
Eigenschaften der positiven tern\"aren quadratischen Formen von Ludwig
August Seber}, G{\"o}ttingische gelehrte Anzeigen, Juli 9 [J. Reine
Angew. Math. 20 (1840), 312-320 = Werke, Vol. 2 (K{\"o}niglische
Gesellschaft der Wissenschaften, G{\"o}ttingen 1876), 188-196].

\bibitem{JK} Kepler, Johannes, \emph{On the six-cornered snowflake},
Oxford Clarendon Press, Oxford 1966 (foreword by L.L.~White) (English
translation of Kepler's Latin essay printed in 1601).

\bibitem{LA} Lagrange, Joseph Louis, \emph{Recherches d'arithmetique},
Nouv. Mem. Acad. Roy. Sc. Belle Letteres, Berlin 1773,  265-312 =
Oeuvres III, 693-758.


\end{thebibliography}
\end{document}